\documentclass[3p,12pt]{elsarticle}
\newcommand{\mysize}{0.7}

\usepackage{graphicx}
\usepackage{amssymb}
\usepackage{amsmath}
\usepackage[colorinlistoftodos]{todonotes}
\usepackage{threeparttable}
\usepackage{color,soul}

\usepackage{lineno}

\begin{document}

\begin{frontmatter}


\title{Optimisation of Wastewater Treatment Strategies in Eco-Industrial Parks: Technology, Location and Transport}



\author[1]{Edward O'Dwyer}
\author[2,3]{Kehua Chen}
\author[2]{Hongcheng Wang}
\author[2]{Aijie Wang}
\author[1]{Nilay Shah}
\author[1]{Miao Guo\corref{cor1}}

\address[1]{Department of Chemical Engineering, Imperial College London, London, SW7 2AZ, UK}
\address[2]{Key Lab of Environmental Biotechnology, Research Center for Eco-Environmental Sciences, Chinese Academy of Sciences, 18 Shuangqing Road, Haidian District, Beijing, 100085, China}
\address[3]{Sino-Danish Center for Education and Research, University of Chinese Academy of Sciences, Beijing, China}
\cortext[cor1]{miao.guo@imperial.ac.uk}

\begin{abstract}
The expanding population and rapid urbanisation, in particular in the Global South, are leading to global challenges on resource supply stress and rising waste generation. A transformation to resource-circular systems and sustainable recovery of carbon-containing and nutrient-rich waste offers a way to tackle such challenges. Eco-industrial parks have the potential to capture symbioses across individual waste producers, leading to more effective waste-recovery schemes. With whole-system design, economically attractive approaches can be achieved, reducing the environmental impacts while increasing the recovery of high-value resources. In this paper, an optimisation framework is developed to enable such design, allowing for wide ranging treatment options to be considered capturing both technological and financial detail. As well as technology selection, the framework also accounts for spatial aspects, with the design of suitable transport networks playing a key role. A range of scenarios are investigated using the network, highlighting the multi-faceted nature of the problem. The need to incorporate the impact of resource recovery at the design stage is shown to be of particular importance. 

\end{abstract}


\end{frontmatter}


\section{Introduction}\label{Intro}

With an increasing population and the acceleration of urbanisation and industrialisation, the appropriation of freshwater resources has increased dramatically in the last few decades \cite{Qu2010,UNESCO2015}, leading in turn to an increase in wastewater. Globally, about \(330km^3/year\) of municipal wastewater is generated \cite{Hoornweg2012} bringing environmental stress, particularly in regions with rapid urbanisation trends. In 2017 in China for example, the discharge of wastewater reached about $7\times10^{6} m^{3}$ \cite{NationalBureauofStatistics2018}. The large quantity of industrial and municipal wastewater has resulted in severe contamination of water bodies (both surface and ground water) \cite{Hu2013}. Responding to this, a number of wastewater treatment (WWT) plants were constructed \cite{NationalBureauofStatistics2014}, with technologies including Anaerobic-Anoxic-Oxic (AAO), oxidation ditch and Sequencing Batch Reactors (SBR) widely adopted \cite{Jin2014a}. In this manner, water quality in China has improved, with national Class 1 - Class 3 water bodies accounting for 67.9\% of the total in 2018.

Despite this progress, the design of treatment strategies that can adopt the most effective new technologies and processes to enhance the removal of contaminants while recovering valuable resources from the wastewater remains an open challenge \cite{Jin2014a}. Waste components of particular importance are phosphorus (P) and nitrogen (N). The N/P discharge to the aqueous environment cause significant eutrophication concerns which can reduce water quality due to the growth of aquatic microbiome e.g. cyanobacteria blooms \cite{Mayer2016}. Furthermore, the global depletion of nutrients, notably nonrenewable P resources, poses a significant challenge to food security due to their essential roles in agricultural sectors. As yet however, resource recovery in the wastewater sector has remained largely under-explored, and represents a research frontier attracting greater attention in recent years \cite{Ashley2011}. Economic viability acts as a significant barrier hindering the implementation of a given resource recovery technology. Thus, it is necessary to develop a system design approach to account for technology integration options as well as the potential cost benefits of waste by-products to explore optimal solutions.

Industrial parks or complexes launched in many countries in Global South offer potential solutions to integrate functional industry networks with eco-efficient design \cite{Farel2016}. Companies and firms can derive economic benefits from land development, construction and shared facilities through industrial parks \cite{Cote1994}. Similarly, the environmental issues faced by the generation of wastewater can be tackled more effectively by exploiting the system symbioses across neighbouring waste generation sites. The presence of different spatially distributed waste streams, with varying compositions and flow characteristics leads to a complex design challenge, particularly if attempting to incorporate resource-circular economic aspects. Topological aspects must be considered for waste transport, while the costs and performances of different technologies need to be included, particularly as new technologies emerge. Progress has been made in this regard, with individual facility optimisation methods developed in \cite{Puchongkawarin2015}. In an eco-industrial park context, optimisation approaches developed to determine optimal resource flows are discussed in \cite{Rubio-Castro2010}, with the incorporation of uncertainty on performance carried out in \cite{Afshari2017}. From a resource recovery perspective, heat recovery networks are considered in \cite{Zhang2016} for example.

Other systems-engineering perspectives can be found for tackling the various design challenges posed. In \cite{Alnouri2016} for example, a network design approach is used to optimise the pipeline connections, allowing for waste stream merging to achieve cost savings in pipework installation. In \cite{Tiu2017} and \cite{Aviso2010} optimisation approaches are presented to enable treatment configurations to be developed that can best exploit water exchange (freshwater and wastewater) between different plants to minimise economic costs and environmental impacts. A need exists to expand this to enable the recovery of a wider range of resources to be accounted for, thus providing a means for emerging recovery technologies to be considered in design. As the literature suggests, the impact of the spatial locations of WWT facilities and the transport networks required to facilitate a given scheme must be included as well as the technology types and configurations. By considering the potential value of resource recovery along with the network design costs and treatment environmental impact, the overall design approach may be significantly altered. Despite research advances on mathematical optimisation in wastewater network design, no research has been published to bridge two-level design challenges i.e. the waste transport (pipeline) and resource recovery technologies.

In this paper, an optimisation framework for design of wastewater treatment schemes in eco-industrial parks is developed. The derived approach allows for any number of treatment technology options to be considered in terms of contaminant removal performance, operational limitations and costs. The technology selection process can allow for wastestreams at different locations in the park to be combined via a waste transport network, the optimal routing for which is also selected within the optimisation approach, the design of which depends on available materials, techniques and topological information. In this way, centralised, decentralised or spatially distributed treatment schemes can be designed, enabling an improved overall performance. The flexible nature of the framework allows for constraints on waste discharge levels for different contaminants, as well as providing the possibility to apply user-defined penalties or taxes on waste discharge in addition to the other system capital and operational costs. A further advantage of the proposed approach is the potential to consider the financial implications of resource recovery, and how the financial performance of a design can be impacted by remuneration through the sale of generated resources. Such a feature can encourage the transition towards a more resource-circular model for wastewater treatment.

In Section \ref{S:2}, the technological considerations for waste treatment in eco-industrial parks are introduced, leading to the description of the optimisation framework (including constraints and objectives) in Section \ref{S:3}. In Section \ref{S:4} a design example and a set of case studies are introduced to illustrate the application of the optimisation framework in a range of different scenarios.

\section{Wastewater Treatment and Industrial Park Design Assumptions}\label{S:2}
\subsection{Current status of resource recovery from wastewater}
Waste stream composition is central to the design of a treatment facility and the selection of suitable technologies. The discharge of different constituent components can impact the wider environment in a variety of ways, while different removal and recovery technologies have varying strengths and weaknesses in targeting specific components. Furthermore, when considering waste as a resource (as alluded to previously), the potential exists to unlock an otherwise untapped revenue stream, however this value is highly dependent on the available market for the specific component in question. Biofuels for example, have a well-established market, whereas the potential for a component such as phosphorus may be less clear \cite{Duan2019}. The costs for discharging an unrecovered resource to the environment (environmental, financial or otherwise) may also prove critical. As different resources require different technologies and the process for extracting one resource may preclude the recovery of another, careful consideration must be made of the various impacts. While a wide range of components could be considered, in this paper we narrow the focus to some of the most environmentally pressing and technologically feasible.
\begin{enumerate}
    \item Phosphorus: P is an essential nutrient to biological growth \cite{Jalali2016}, while the discharge of phosphorus can cause eutrophication and the deterioration of water bodies \cite{Lu2016a}. In addition, the natural supplies of phosphate rock are non-renewable and will be exhausted in 30-300 years \cite{Mew2016}. To reduce the depletion of phosphate rocks and prevent eutrophication in natural water bodies, the recovery of phosphorus from wastewater is an effective method \cite{Loganathan2014}. To date, several technologies are mainly applied for phosphorus recovery: chemical precipitation \cite{Jones2015}, adsorption \cite{Loganathan2014}, wet-chemical treatment \cite{Appels2010}, thermochemical treatment \cite{Adam2009}.
    \item Nitrogen: nitrogen is another vital element to all organisms and its discharge can also lead to eutrophication \cite{Huang2017a}. The recovered nitrogen can be used as fertilizer replacement to agricultural production \cite{Kelly2014}. To our knowledge, the recovery of nitrogen can be achieved by bioelectrochemical system \cite{Pandey2016}, ion exchange \cite{Smith2015}, membrane technology \cite{Rahaman2014}, biological technology (e.g. microalgae, \cite{Cai2013}) and struvite crystallization \cite{Huang2016}.
    \item Bio-fuels from sludge recovery: sewage sludge is a by-product of wastewater treatment \cite{Feng2014}. There are heavy metals, pathogens, organic contaminants and enriched biomass in sludge \cite{Dong2013}, which means sludge has complex components but high reused values. At present, there are a few options to recover resources from sludge \cite{Rulkens2008}: anaerobic digestion, production of biofuels, electricity production through microbial fuel cells, incineration for energy recovery, pyrolysis and gasification, utilization for construction materials, supercritical wet oxidation, hydrothermal treatment. Various technologies can acquire different recovered products.
\end{enumerate}	

\subsection{Transport Network Design}
The transport of waste in an industrial park is an important aspect in the design of an appropriate park-wide strategy. Combining wastestreams in centralised plants can only be considered if it is economically viable to construct and operate a transport network between the sites. A highly detailed model of the pipeline construction and associated costs would not be practical for a high-level design tool such as that proposed here however, and as such, certain simplifying assumptions are made. Trench installation (TI) methods are considered using high-density polyethylene pipework (HDPE). The cost of the pipe network installation is then assumed to be a function of the pipe diameter, the elevation change between the start and end of the pipe and the conveyance distance. Using the Chinese Water Supply and Drainage Design Manual \cite{Chinamunicipal2014}, the cost per unit length of different construction methods in discrete depths with specific materials and diameters can be found. Labour and machinery costs are assumed to be constant per unit length.

\section{Optimisation Formulation}\label{S:3}
The complete optimisation formulation is presented in this section including all constraints and objectives. The problem is formulated as a Mixed Integer Linear Program (MILP), with linearisation of some constraints required. Sets and indices are defined in Table \ref{Tab:SetDefine}.

\begin{table}[htpb]
\centering
\begin{tabular}{l|l}
\hline
 & \textbf{Description}\\
\hline\hline
\(j\in J\) & Denotes individual wastewater stream\\
\(p\in P\) & Components associated with wastestream\\
\(r\in R\) & Resources recovered at treatment facility\\
\(x\in X\) & Denotes grid cell\\
\(m\in M\) & Treatment technology or pipeline connector\\
\(JM\subset M\) & Set of connector types\\
\(l\in L\) & Pipeline type\\
\(el\in El\) & Cell elevations present in park\\
\(t\in T\) & Time-steps in operational period\\
\hline
\end{tabular}
\caption{Sets and indices used in the formulation}
\label{Tab:SetDefine}
\end{table}

\subsection{Treatment and technology constraints}

\begin{figure}[htbp]
    \centering
    \includegraphics[width=\mysize\textwidth]{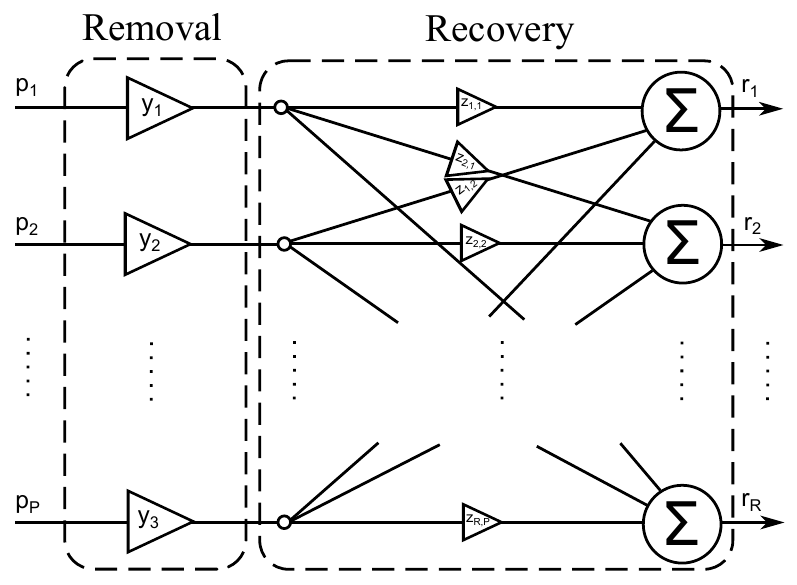}
    \caption{Generic schematic of removal and recovery of resources in a treatment plant}
    \label{fig:Schematic}
\end{figure}

In cell \(x\) for wastestream \(j\), the presence of a treatment plant or pipe connector of type \(m\) is denoted by the binary variable \(\alpha_{j,x,m}\), while the presence of a waste generating facility is given by the binary variable \(G_{j,x}\). Such a facility is assumed to generate waste with component characteristics \(p\) given at time \(t\) by \(gen_{j,x,p,t}\). The quantity per unit flow of each component \(p\) in stream \(j\) (in the case of chemical constituents, this can represent the concentration for example) is given by \(C_{j,p}\), with the flowrate of the stream from cell \(x\) to cell \(x^{'}\) at time \(t\) given as \(fl_{j,x,x^{'},t}\). With this notation, the change in each component resulting from treatment in plant type \(m\) is given by:

\begin{gather}
    \Delta P_{j,m,x,p,t}=\sum_{x'}fl_{j,x',x,t}C_{j,p}\alpha_{j,m,x}y_{m,p},\\\nonumber \forall j\in J,\forall m\in M,\forall x\in X, \forall p\in P, \forall t\in T
\end{gather}

where \(y_{m,p}\) can be seen as the efficiency with which a plant \(m\) removes component \(p\). In some plants, the recovery of useful resources may be possible, with the quantity of recovered resource \(r\) denoted \(Rec_{m,x,r,t}\). The recovery rate of a resource can be characterised for a treatment plant by a recovery efficiency, whereby a certain proportion of a removed inflow component is recovered as a useful resource. More generally, to allow for the generation of resources not present in the inflow, the resources are assumed to be generated at a rate that is proportional to some linear combination of the removed components. This is captured by the following constraint:

\begin{gather}
    Rec_{j,m,x,r,t}=\sum_{p\in P}z_{m,r,p}\Delta P_{j,m,x,p,t},\\\nonumber \forall j\in J,\forall m\in M,\forall x\in X, \forall r\in R, \forall t\in T
\end{gather}

where \(z_{m,r,p}\) maps the relationship between removed component \(p\) and recovered resource \(r\) for treatment plant type \(m\). A schematic of a generic treatment plant type is shown in Fig.\ref{fig:Schematic}. In the diagram, the plant is divided into a removal and a recovery section. In the former, the quantity of each component that is removed from the flow and thus available for recovery is calculated and passed to the recovery section where generation of resources is carried out. The recovered resources may be the of the same type as the removed components (e.g. phosphorus (P) could be viewed as a waste component and a useful resource) or different (e.g. Chemical Oxygen Demand (COD) could be a characteristic of the waste and can be converted to methane ($CH_{4}$) for use as a biofuel). 

The total flow from all wastestreams to a treatment plant \(m\) in cell \(x\) is bounded by the limit \(F^{max}\) (the maximum flow capacity) as follows:

\begin{gather}
    \sum_{j\in J}\sum_{x'\in X}fl_{j,x',x,t}\alpha_{j,m,x}\leq F^{max}_{m},\quad\forall m\in M,\forall x\in X,\forall t\in T
\end{gather}

The presence of any treatment plant for wastestream \(j\) in cell \(x\) is indicated by the binary \(\upsilon_{j,x}\), determined by the following constraints ($\forall j\in J,\forall x\in X$):

\begin{gather}
    \mathcal{M}\upsilon_{j,x}\geq\sum_{m\in M\setminus JM}\alpha_{j,x,m}\\
    \mathcal{M}\left(1-\upsilon_{j,x}\right)>-\sum_{m\in M\setminus JM}\alpha_{j,x,m}
\end{gather}

where \(\mathcal{M}\) represents an arbitrary large number such that \(\mathcal{M}>>|X|\).
Similarly, the presence of a particular treatment plant \(m\) in cell \(x\) for any wastestream is indicated by the binary \(\omega_{x,m}\), given by the constraints:

\begin{gather}
    \mathcal{M}\omega_{x,m}\geq\sum_{j\in J}\alpha_{j,x,m}\\
    \mathcal{M}\left(1-\omega_{x,m}\right)>-\sum_{j\in J}\alpha_{j,x,m} 
\end{gather}

The total quantity of recovered resources \(r\) from the park at time \(t\) is given as:

\begin{gather}
    rs_{r,t}=\sum_{j\in J}\sum_{x\in X}\sum_{m\in M}Rec_{j,m,x,r,t},\quad\forall r\in R,\forall t\in T
\end{gather}

Similarly, the total quantity of discharged component \(p\) from the park is:

\begin{gather}
    dis_{j,p,t}=\sum_{x\in X}\left(gen_{j,x,p,t}-\sum_{m\in M}\Delta P_{j,m,x,p,t}\right),\quad\forall j\in J,\forall p\in P,\forall t\in T
\end{gather}

Explicit park discharge restrictions can be applied to each component \(p\):

\begin{gather}
    \sum_{j\in J}dis_{j,p,t}-P^{max}_{p}\leq 0,\quad\forall p\in P,\forall t\in T
\end{gather}

\subsection{Transport constraints}
\begin{figure}[htbp]
    \centering
    \includegraphics[width=\mysize\textwidth]{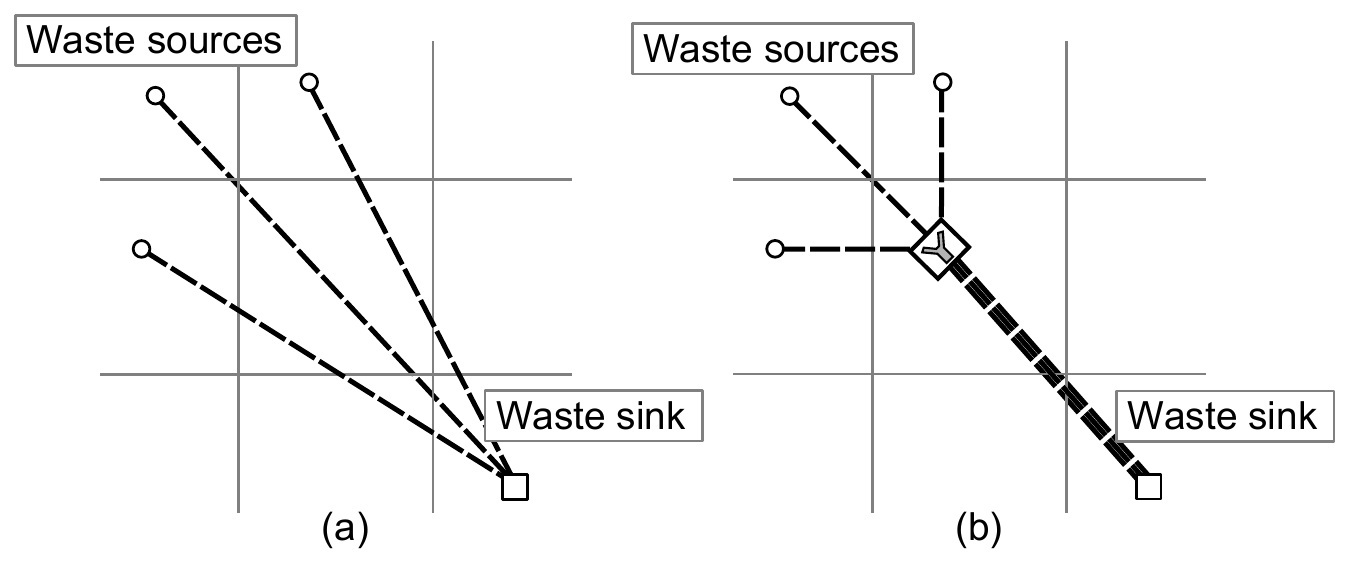}
    \caption{Transport pathway design options - waste streams can use a common tranpsort pathway to reduce construction costs}
    \label{fig:Joints}
\end{figure}

In this section, the constraints associated with the transport of waste through pipeline installations in the park are presented. To reduce pipe trench installation costs, the formulation allows for joints to be included in the network for flows to be transported in the same trench. This is shown in Fig.\ref{fig:Joints}, where transport of waste from a source to a sink is illustrated. In Fig.\ref{fig:Joints}(a) the waste is directly transported, whereas in  Fig.\ref{fig:Joints}(b), a joint is used to reduce the total trench installation requirements. 

For each pipeline type \(l\), the flowrate is bounded by the limit \(Fl^{max}_{l}\) calculated as:

\begin{gather}
    Fl^{max}_{l}=vel_{l}Area_{l},\quad\forall l\in L
\end{gather}

where \(vel_{l}\) is the design velocity and \(Area_{l}\) is the cross-sectional area of the pipe. The selection of a particular pipe type is indicated by the binary variable \(\delta_{l}\). The total flow of waste out of cell \(x\) is then limited by the constraint:

\begin{gather}
    \sum_{j\in J}\sum_{x'\in X}fl_{j,x,x',t}\leq\sum_{l\in L}Fl^{max}_{l}\delta_{l},\quad\forall x\in X,\forall t\in T\\
    \sum_{l\in L}\delta_{l}=1
\end{gather}

Waste can flow from a generation site or from a pipeline connection. As such, the total flow out of any cell must be bounded by the total flow into a connector in the cell plus the generated waste in that cell. This leads to the following constraint which describes this limiting condition for waste component \(p\) at time \(t\):

\begin{gather}
    \sum_{x'\in X}fl_{j,x,x',t}C_{j,p}\leq\sum_{x'\in X}\sum_{m\in JM}\alpha_{j,x,m}fl_{j,x',x,t}C_{j,p}+gen_{j,x,p,t},\\\nonumber\quad\forall j\in J,\forall x\in X,\forall p\in P,\forall t\in T
\end{gather}

Similarly, the flow into a connector must also flow out again, as enforced by the following constraint:

\begin{gather}
    \sum_{m\in JM}\sum_{x'\in X}fl_{j,x',x,t}\alpha_{j,x,m}\leq\sum_{x'\in X}fl_{j,x,x',t},\quad\forall j\in J,\forall x\in X,\forall t\in T
\end{gather}

The presence of a connector at any point in a wastestream \(j\) is indicated by the binary variable \(\epsilon_{j}\), determined $\forall j\in J$ by:

\begin{gather}
    \mathcal{M}\epsilon_{j}\geq\sum_{m\in JM}\sum_{x\in X}\alpha_{j,x,m}\\
    \mathcal{M}\left(1-\epsilon_{j}\right)>-\sum_{m\in JM}\sum_{x\in X}\alpha_{j,x,m}
\end{gather}

The presence of a waste generation site in a cell \(x\) is indicated by the pre-defined binary parameter \(G_{j,x}\), while the presence of transport pathway between a generation cell \(x\) (\(\forall x\in X\)) and cell \(x'\) (\(\forall x'\in X\)) with a pipeline connection is then indicated by \(\varphi_{x,x'}\), defined by the following constraints:

\begin{gather}
    \mathcal{M}\varphi_{x,x'}\geq\sum_{j\in J}\sum_{m\in JM}\alpha_{j,x',m}G_{j,x}\\
    \mathcal{M}\left(1-\varphi_{x,x'}\right)>-\sum_{j\in J}\sum_{m\in JM}\alpha_{j,x',m}G_{j,x}\\
    \varphi_{x,x'}\leq \sum_{j\in J}\sum_{m\in JM}\alpha_{j,x',m}\\
    \varphi_{x,x'}\leq \sum_{j\in J}G_{j,x}
\end{gather}

A transport pathway between a cell with a pipeline connection \(x\) (\(\forall x\in X\)) and a cell \(x'\) (\(\forall x\in X\)) with a treatment plant \(m'\) is indicated by the binary variable \(\Pi_{x,x'}\). This is set by the following constraints (\(\forall j\in J\)), where \(\kappa_{j,x,x'}\) is a binary variable used as a dummy intermediate variable.

\begin{gather}
    \mathcal{M}\kappa_{j,x,x'}\geq\sum_{m\in JM}\sum_{m'\in M\setminus JM}\alpha_{j,x,m}\alpha_{j,x',m'}\\
    \mathcal{M}\left(1-\kappa_{j,x,x'}\right)>-\sum_{m\in JM}\sum_{m'\in M\setminus JM}\alpha_{j,x,m}\alpha_{j,x',m'}\\
    \mathcal{M}\Pi_{x,x'}\geq\sum_{j\in J}\kappa_{j,x,x'}\\
    \mathcal{M}\left(1-\Pi_{x,x'}\right)>-\sum_{j\in J}\kappa_{j,x,x'}
\end{gather}

Finally, the presence of a transport pathway directly between a generation cell \(x\) (\(\forall x\in X\)) and a treatment cell \(x'\) (\(\forall x\in X\)) without passing through a connection is indicated by the binary \(\gamma_{x,x'}\), where:

\begin{gather}
    \mathcal{M}\gamma_{x,x'}\geq\sum_{j\in J}G_{j,x}\upsilon_{j,x'}\left(1-\epsilon_{j}\right)\\
    \mathcal{M}\left(1-\gamma_{x,x'}\right)>-\sum_{j\in J}G_{j,x}\upsilon_{j,x'}\left(1-\epsilon_{j}\right)
\end{gather}

Bringing this together, a transport pathway of any type between any cell \(x\) and \(x’\) is denoted \(\partial_{x,x'}\), as given by:

\begin{gather}
    \mathcal{M}\partial_{x,x'}\geq\left(\varphi_{x,x'}+\gamma_{x,x'}+\Pi_{x,x'}\right)
\end{gather}

\subsection{Objectives}
The overall objective for the problem is to minimise the total cost which can be taken as a combination of the transport costs (\(Cost^l\)), the treatment costs (\(Cost^t\)), the resource discharge penalties (\(Cost^d\)) and the resource sales (\(Sell^r\)) as follows:

\begin{gather}
    J = Cost^l+Cost^t+Cost^{d}-Sell^{r}
\end{gather}

To determine the transport costs, a parameter (\(Pipe^{C}_{x,x',l}\)) describing the cost of a pipeline of type \(l\) installed between cells \(x\) (\(\forall x\in X\)) and \(x'\) (\(\forall x'\in X\)) can be found as:

\begin{gather}
    Pipe^{C}_{x,x',l} = \sum_{el\in El}\mu_{x,x',el}\left(Cap^{pipe}_{l,el}+Op^{l}_{el,l}\right)\ell_{x,x'},\quad\forall l\in L
\end{gather}

where \(Cap^{pipe}_{l,el}\) is the installation cost per length of pipe of type \(l\), \(Op^{l}_{el,l}\) is the operational cost of pumping waste, \(\ell_{x,x'}\) is the distance between \(x\) and \(x'\) and \(\mu_{x,x',el}\) indicates the elevation change, given as (\(\forall x\in X, \forall x'\in X,\forall el\in El\)):

\begin{gather}
    \mu_{x,x',el}=\begin{cases}
    1 & \Delta_{x,x'} = el\\0 & otherwise
    \end{cases}
\end{gather}

where \(\Delta_{x,x'}\) is the elevation difference between \(x\) and \(x'\).
From this, the installation cost of the optimal transport route is given by:

\begin{gather}
    Cost^l=\sum_{l\in L}\sum_{x\in X}\sum_{x'\in X}Pipe^{C}_{x,x',l}\delta_{l}\partial_{x,x'}
\end{gather}

The operational costs of waste treatment are given (\(\forall x\in X,\forall m\in M\)) by:

\begin{gather}
    CostOp^{t}_{x,m}=\sum_{j\in J}\sum_{x'\in X}\sum_{t\in T}\alpha_{j,m,x}fl_{j,x',x,t}Op^{t}_{m}
\end{gather}

where \(Op^{t}_{m}\) is the operational cost per unit of treated waste for plant type \(m\) in cell \(x\). The capital cost of such a treatment plant is given (\(\forall x\in X,\forall m\in M\)) by:

\begin{gather}
    CostCap^{t}_{x,m}=Cap^{t}_{m}\omega_{x,m}
\end{gather}

where \(Cap^{t}_{m}\) is the cost of installing plant type \(m\). The total treatment cost is then:

\begin{gather}
    Cost^{t}=\sum_{x\in X}\sum_{m\in M}\left(CostCap^{t}_{x,m}+CostOp^{t}_{x,m}\right)
\end{gather}

To monetise the environmental impacts of excessive waste discharge, additional discharge penalty costs can be applied to waste component \(p\), using the parameter \(Pen_{p}\), leading to a total penalty cost of:

\begin{gather}
    Cost^{fin}=\sum_{j\in J}\sum_{p\in P}\sum_{t\in T}\left(Pen_{p}dis_{j,p,t}\right)
\end{gather}

Conversely, income can be generated by applying a selling price, \(price_{r}\), to recovered resource \(r\). The total income generated in the park is then given as:

\begin{gather}
    Sell^{r}=\sum_{r\in R}\sum_{t\in T}\left(rs_{r,t}price_{r}\right)
\end{gather}

\section{Design Case}\label{S:4}
\subsection{Description of illustrative park}
To demonstrate the effective application of the optimised treatment facility design, a set of hypothetical case studies are presented in this section. The purpose is to illustrate how different factors influence the optimal design of a treatment scheme and as such, how a comprehensive approach that considers transport, topological features, circular economic aspects and technological performance characteristics is required to enable effective design. The ability of the proposed methodology to incorporate these factors is shown. In all cases, the optimisation problem was solved using the CPLEX solver (version 25.0.3) implemented using the GAMS software.

In the hypothetical design case, a set of five different waste streams generated from diverse industrial facilities are considered, at different locations within the same industrial park. In the design cases carried out here, waste is characterised only by its COD, Total Nitrogen (TN) content and Total Phosphorus (TP) content, however, the modelling framework has no restrictions on what parameters can be considered. The hypothetical design case is based on the operational data collected from an industrial park located in East China. As shown in Table \ref{Tab:Comps}, the compositions of the wastewater streams vary significantly with different industries. The selected industrial wastewater ranges from low COD effluent discharged from printing and mixed industry sources to COD-rich effluent from fermentation and food processing industries e.g. corn processing and pharmaceutical industries. 


\begin{table}[!htpb]
\centering
\begin{threeparttable}

\begin{tabular}{l c c c c c}
\hline
 & \textbf{Source} & \textbf{Flowrate} & \textbf{COD} & \textbf{TN} & \textbf{TP} \\
 & & \(m^3/d\) & \(mg/L\) & \(mg/L\) & \(mg/L\)\\
\hline\hline
Stream A & Domestic and industrial mixed & 10000 & 713 & 86.3 & 0.4 \\
Stream B & Printing and dying industry & 10000 & 400 & 40 & 7 \\
Stream C & Corn processing industry & 8000 & 1500 & 100 & 80 \\
Stream D & Pharmaceutical industry & 4000 & 2030 & 126.3 & 1.74 \\
Stream E & Chinese medicine industry & 3000 & 15000 & 420 & 245 \\
\hline
\end{tabular}
\caption{Waste stream compositions from 5 industrial sources}
\label{Tab:Comps}
    
 \end{threeparttable}
\end{table}

\begin{figure}[htbp]
    \centering
    \includegraphics[width=\mysize\textwidth]{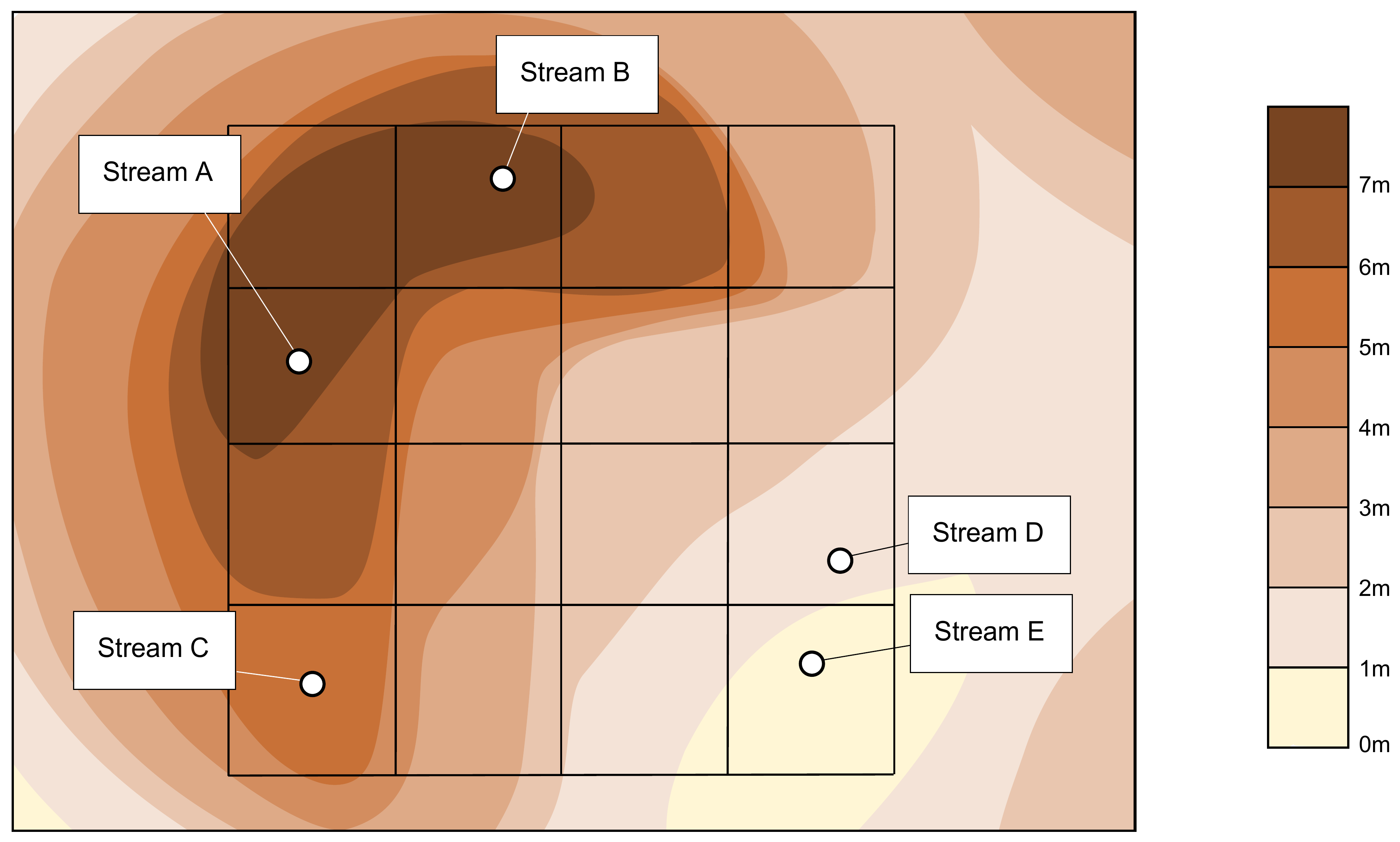}
    \caption{Topology of industrial park overlaid by 500m grid}
    \label{fig:Contours}
\end{figure}

The different waste producers are located within an area of 2\(km^2\) on a landscape of varying elevation. To formulate the MILP problem, the park is split into a $4\times4$ mesh of cells, with the length of each vertex set at \(500m\). It should be noted that the number of cells, as well as the cell size, shape and distribution can be chosen arbitrarily (only the distance and elevation change between cells is required in the formulation). The park topology and the locations of the sources of the waste streams are shown in Fig. \ref{fig:Contours}, overlaid by outlines of the chosen cells.

\subsection{Technology Options}


To demonstrate the model functionality, we developed a hypothetical study considering technologies with different scales and performance characteristics to remove and recover COD, TN and TP. As given in Table \ref{Tab:Techs}, different technology options were modelled i.e. aerobic WWT facility, A2/O WWT (anaerobic-anoxic-aerobic), integrated anaerobic digestion and microalgae cultivation. Technologies A-D represent larger-scale treatment facilities favouring more centralised schemes, while E-H represent smaller distributed versions of the same technologies. Technology options A and E are based on a WWT plant adopting an anaerobic hydrolysis acidification-circulating activated sludge system and denitrification filter-fiber turntable filter-disinfection. Performance characteristics were derived from a domestic WWT facility located in Zhejiang Dongyang in China (based on unpublished data), treating streams of 8\% industrial wastewater and 92\% domestic wastewater. Technologies B and F are based on integrated anaerobic digestion and algae cultivation, with removal and recovery efficiencies taken from \cite{Arias2018}. The technology options C and G were adapted from domestic WWT plant data, which is a circulating activated sludge system located in Jiansu Taizhou, China. The technologies D and H are based on the A2/O and Membrane Bioreactor (MBR) system located in Beijing China where the removal efficiencies were derived from the whole WWT processes. The total costs for each technology option were estimated based on the assumption that 50\%\ costs were caused by operational and capital costs respectively, where the operational costs were derived from industrial on-site data. 


\begin{table}[htpb]
\centering
\begin{threeparttable}
\begin{tabular}{c| c| c c c| c c c}
\hline
 & & \multicolumn{3}{c|}{\textbf{Removal \(\eta\)}} & \multicolumn{3}{c}{\textbf{Recovery \(\eta\)}}\\\hline& \textbf{Cap} & \textbf{COD} & \textbf{TN} & \textbf{TP} & \textbf{\(CH_{4}\)} & \textbf{N} & \textbf{P}\\
 & \(m^3/d\) & \(\%\) & \(\%\) & \(\%\) & \(L_{rec}/g_{COD_{removed}}\) & \(\%\) & \(\%\) \\
\hline\hline
A\tnote{1} & 40000 & 93 & 63 & 87 & 0.375 & 8 & 87 \\
B\tnote{2} & 40000 & 70 & 100 & 100 & 0.596 & 0 & 0 \\
C\tnote{3} & 40000 & 88 & 84 & 30 & 0 & 0 & 0 \\
D\tnote{4} & 40000 & 93 & 75 & 97 & 0 & 0 & 0 \\
E\tnote{5} & 10000 & 93 & 63 & 87 & 0.375 & 8 & 87 \\
F\tnote{6} & 10000 & 70 & 100 & 100 & 0.596 & 0 & 0 \\
G\tnote{7} & 10000 & 88 & 84 & 30 & 0 & 0 & 0 \\
H\tnote{8} & 10000 & 93 & 75 & 97 & 0 & 0 & 0 \\
\hline
\end{tabular}
\caption{Removal and recovery efficiency of technology options considered}
\label{Tab:Techs}

\begin{tablenotes}
\small
  \item[1,5] {Data derived from experiments at Chinese Academy of Sciences}
   \item[2,6]{Facility assumed to integrate WWT and microalgae cultivation}
   \item[3,7]{Aerobic WWT plants }
   \item[4,8]{Anaerobic-anoxic-aerobic }
  \end{tablenotes}

 \end{threeparttable}
\end{table}
\color{black}
\subsection{Pipework Costs}
The transport of waste can be carried out using high-density polyethylene (HDPE) pipework with possible diameters ranging from \(300mm-600mm\). To capture a more complete view of the installation cost of the pipework, the trenches into which the pipes are laid must be taken into account. To maintain suitable gradients, deeper trenches are required for larger differences in elevation between source and sink. Table \ref{Tab:Pipes} summarises the cost of the different pipework options (taken from \cite{Chinamunicipal2014}) across different elevation changes.
\begin{table}[htpb]
\centering
\begin{tabular}{c| c c c c}
\hline
 \(\Delta\) Elev & \multicolumn{4}{c}{Cost per \(100m\) (\pounds)} \\
 \(m\) & \(\Phi=0.3m\) & \(\Phi=0.4m\) & \(\Phi=0.5m\) & \(\Phi=0.6m\) \\
\hline\hline
0 & 275 & 465 & 809 & 1111 \\
1.5 & 3765 & 4145 & 4830 & 5434 \\
2.5 & 4944 & 5324 & 6009 & 6612 \\
3.5 & 5478 & 5857 & 6541 & 7144 \\
4.5 & 29248 & 29626 & 30308 & 30910 \\
5.5 & 46111 & 46487 & 47166 & 47764 \\
6.5 & 114498 & 114873 & 115550 & 116147 \\
7.5 & 116165 & 116541 & 117218 & 117814 \\
\hline
\end{tabular}
\caption{Cost per 100m length of pipework installation for different elevation changes}
\label{Tab:Pipes}
\end{table}

The flow through the pipes is limited by a flow velocity constraint, in this case set at \(2m/s\) as well as a pipe capacity constraint, set here at \(80\%\).

\subsection{Case study scenario selection}
A set of case studies are presented in this section to illustrate the capability of the optimisation framework to incorporate different objectives in determining optimal design schemes. The potential financial benefits available through the implementation of suitable transport networks and centralised treatment plants are first demonstrated. Following this, the ability to impose strict limits on the discharge of different resources is shown. Finally, the design impact of the transition to a waste circular framework in which valuable resources can be recovered is highlighted. In all cases, 10 years of operation are assumed (an arbitrary value that can be chosen by the user).

\subsubsection{From decentralised to distributed}\label{Sec:Distributed}
\begin{figure}[htbp]
    \centering
    \includegraphics[width=\mysize\textwidth]{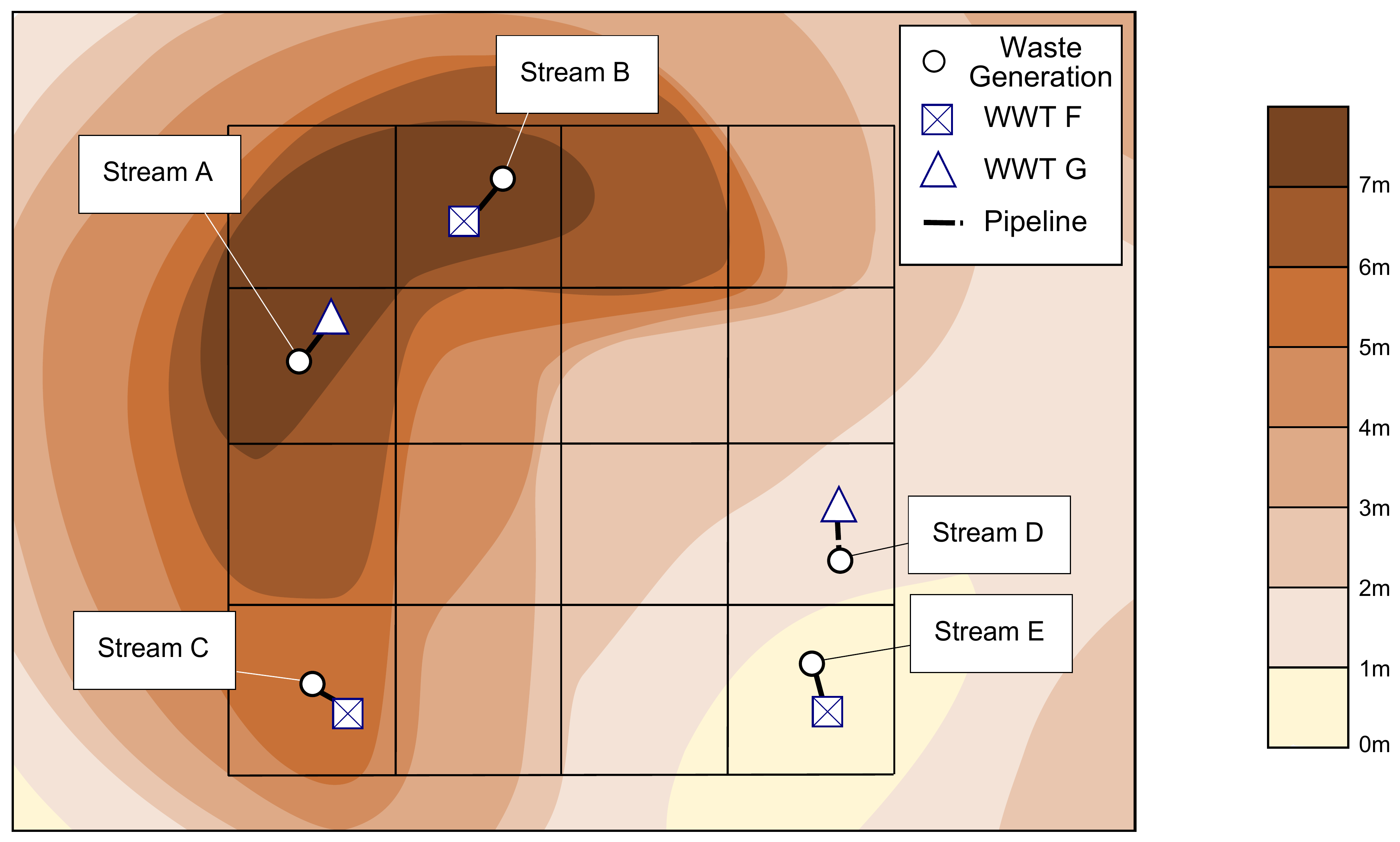}
    \caption{Wastewater treatment scheme without transport}
    \label{fig:ContoursDecentralised}
\end{figure}

 The optimisation framework is first used to determine the optimal technology selection for a decentralised scenario in which no transport is possible. To incentivise the removal of N and P from discharged wastewater streams in the park, a penalty of \(\pounds0.8\) is applied to each \(kg\) of N or P not removed from each stream. The resulting design scheme is shown in Fig.\ref{fig:ContoursDecentralised} where it can be seen that three F-type plants and two G-type plants are required, each operating below maximum capacity.

\begin{figure}[htbp]
    \centering
    \includegraphics[width=\mysize\textwidth]{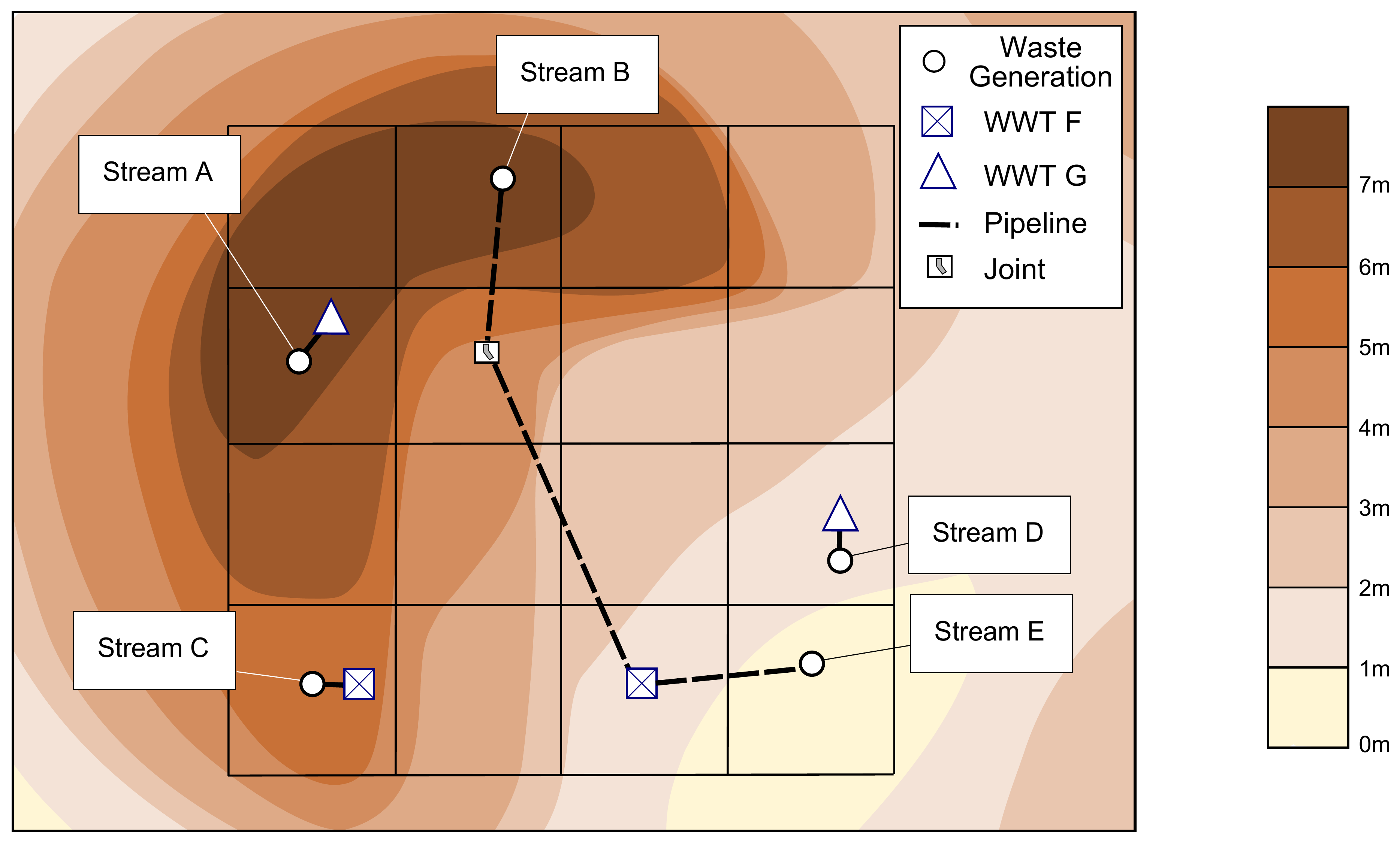}
    \caption{Optimal treatment scheme layout with transport enabled}
    \label{fig:ContoursDisCost}
\end{figure}

\begin{figure}[htbp]
    \centering
    \includegraphics[width=\mysize\textwidth]{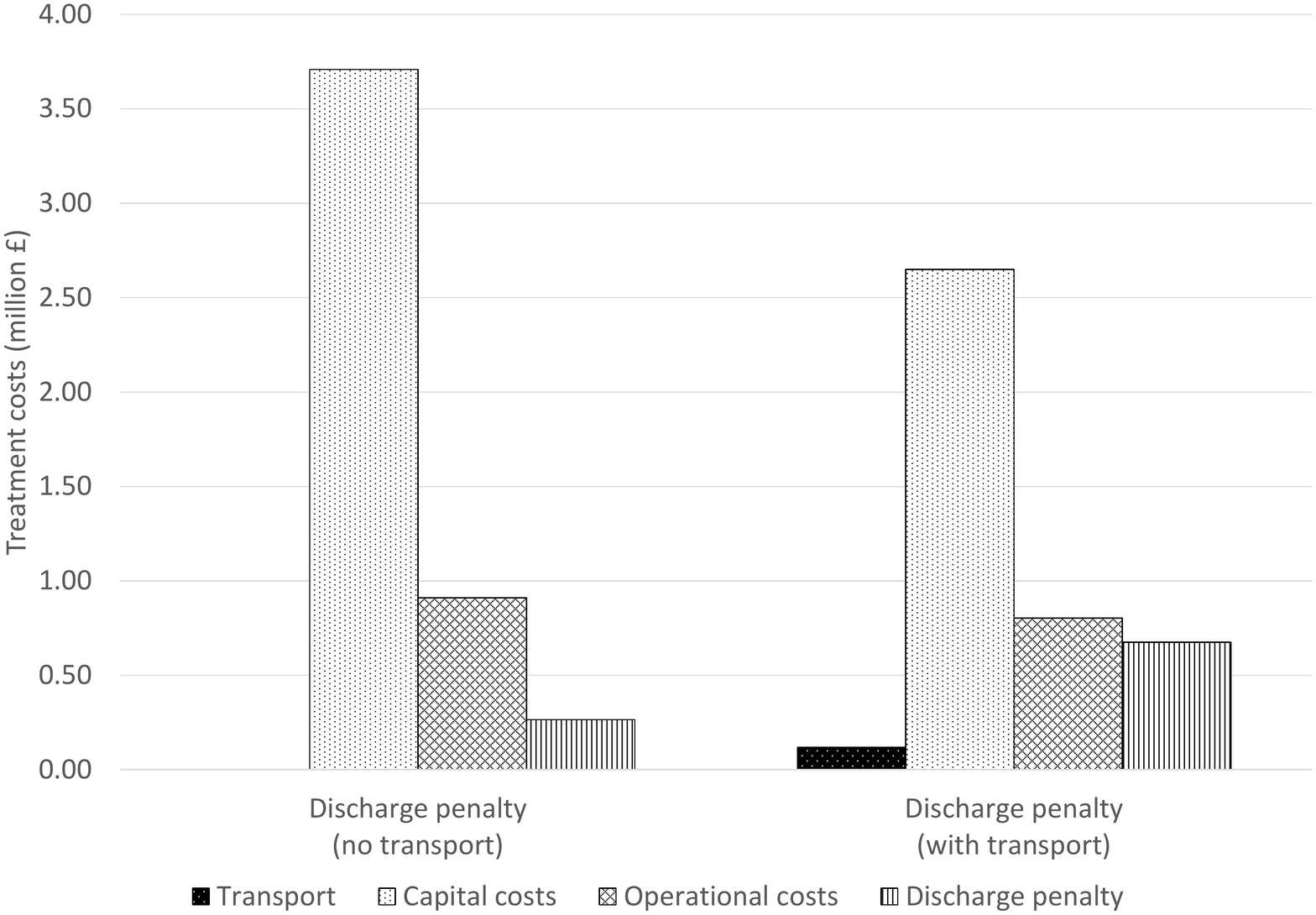}
    \caption{Financial breakdown of treatment costs for eco-industrial parks with and without transport networks}
    \label{fig:BarDecentVsDistFin}
\end{figure}

By allowing transport of waste across the park, the plants can operate more closely to their operational capacities, thereby achieving a more cost-effective outcome. The optimal design scheme is shown in Fig.\ref{fig:ContoursDisCost}. In this case, Stream B and Stream E are transported to a common waste treatment centre. The high volume and low COD and TN concentration characteristics of Stream B combine with the low volume, high concentration characteristics of Stream E to ensure the best utilisation of the treatment plant. The remaining waste generators carry out treatment on site. The total cost for construction and 10 years of operation of the wastewater treatment facility was found to be 13\% lower with transport included, with a breakdown provided in Fig.\ref{fig:BarDecentVsDistFin}.

\subsubsection{Environmental objectives and the economic impact}
Apart from the financial aspects of the treatment schemes, the environmental impacts must also be considered. In the previous examples, a penalty was applied to discourage the discharge of excess pollutants to the environment. This penalty could be adjusted to redress the balance between environmental and financial objectives, however, the formulation presented here also allows for the explicit inclusion of strict limits on the discharge of the different pollutants. To demonstrate this, two new scenarios are developed here within which the discharge penalties are removed and strict thresholds on the total discharge of N and P from the park are applied. For Scenario 1, typical standards are applied \cite{ISO15681,En2008}, while for Scenario 2, a more environmentally focussed design is desired, and enforced by reducing these discharge limits by a factor of 10. The values are shown in Table \ref{Tab:DischargeLimits}.

\begin{table}[htpb]
\centering
\begin{tabular}{c| c c}
\hline
& N discharge limit\tnote{1} & P discharge limit\tnote{2} \\
 & \(kg/L\) & \(kg/L\) \\
\hline\hline
Scenario 1 & \(1.5\times10^{-5}\) & \(1\times10^{-6}\) \\
Scenario 2 & \(1.5\times10^{-6}\) & \(1\times10^{-7}\) \\
\hline
\end{tabular}
\caption{Limits applied to discharge of N and P from the park}
\label{Tab:DischargeLimits}
  
  
\end{table}

\begin{figure}[htbp]
    \centering
    \includegraphics[width=\mysize\textwidth]{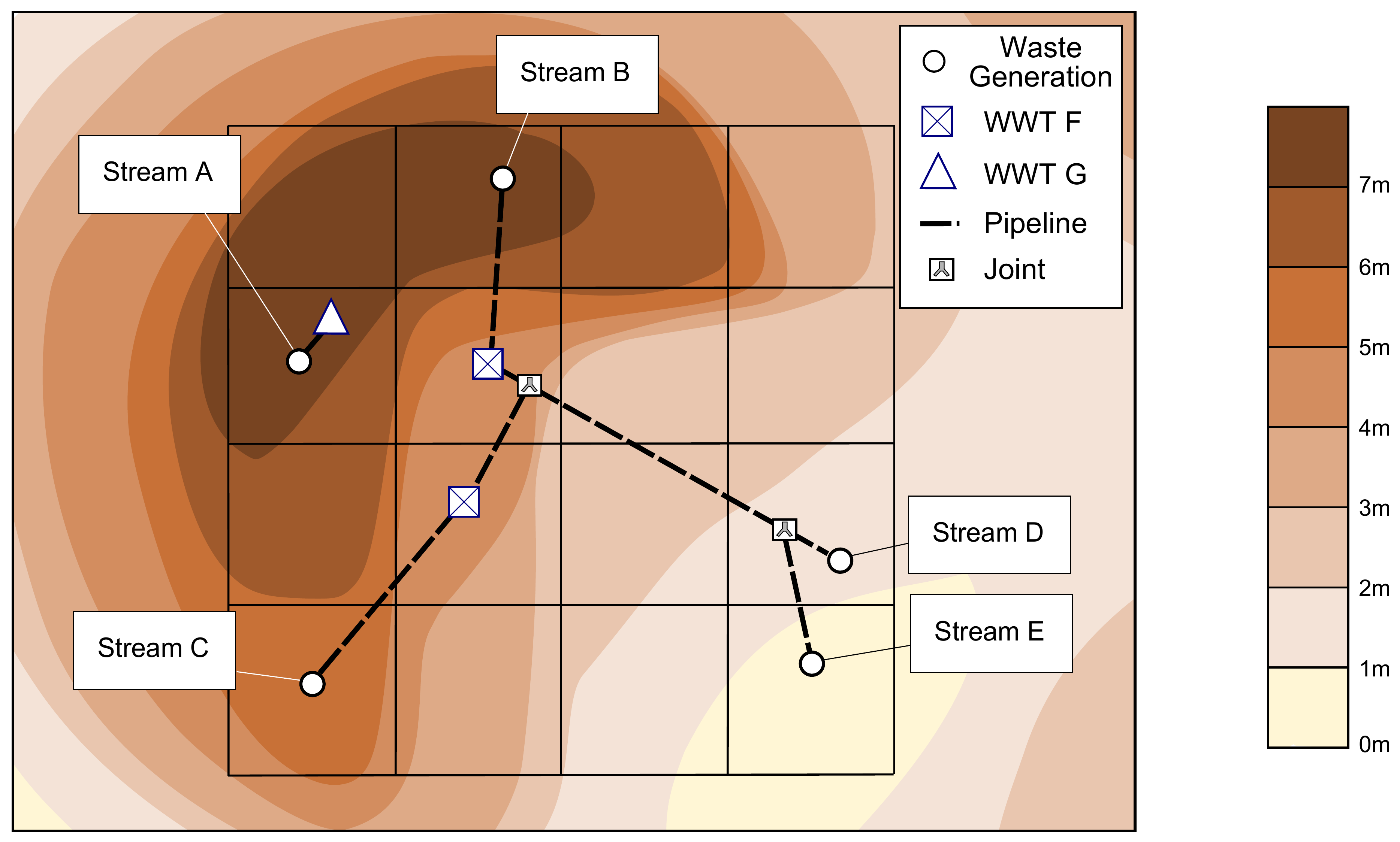}
    \caption{Optimal treatment scheme layout with transport network}
    \label{fig:ContoursNoDisCost}
\end{figure}
\begin{figure}[htbp]
    \centering
    \includegraphics[width=\mysize\textwidth]{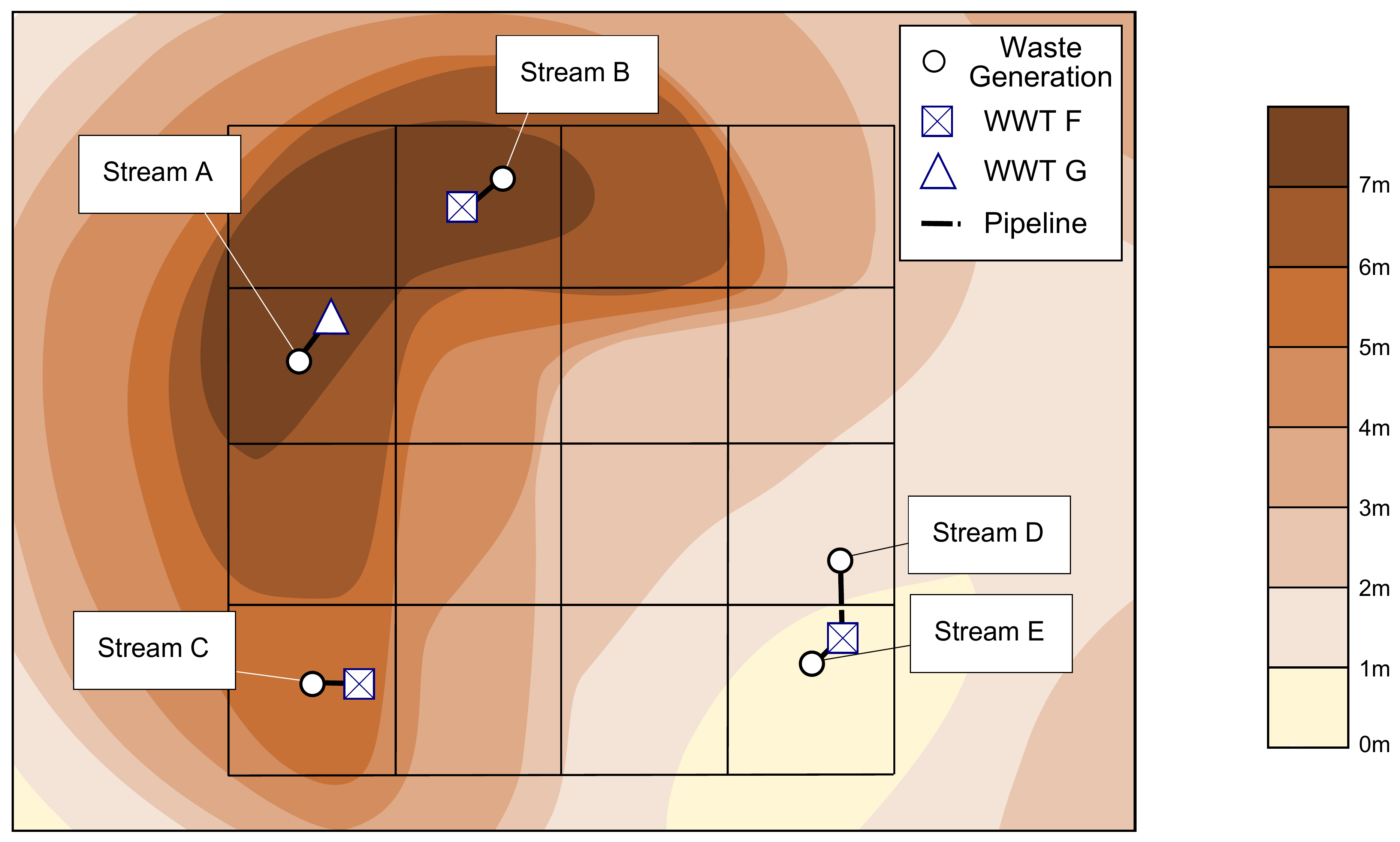}
    \caption{Optimal treatment scheme layout with transport network and strict discharge limits}
    \label{fig:ContoursNoDisCostTenth}
\end{figure}

The resulting schemes for Scenario 1 and Scenario 2 are shown in Fig. \ref{fig:ContoursNoDisCost} and Fig. \ref{fig:ContoursNoDisCostTenth} respectively. These schemes are quite different - when stricter limits are applied, a more decentralised scheme is chosen, with most streams carrying out treatment at the point of waste generation, ensuring that adequate treatment capacity is available to treat all waste. Only Stream D and Stream E combine - their lower volumes enable a single plant to be used without exceeding the plants capacity limit. In Scenario 1, waste is transported towards the centre of the park - less treatment plant capacity is installed, leading to reduced costs, but larger quantities of untreated waste.

\begin{figure}[htbp]
    \centering
    \includegraphics[width=\mysize\textwidth]{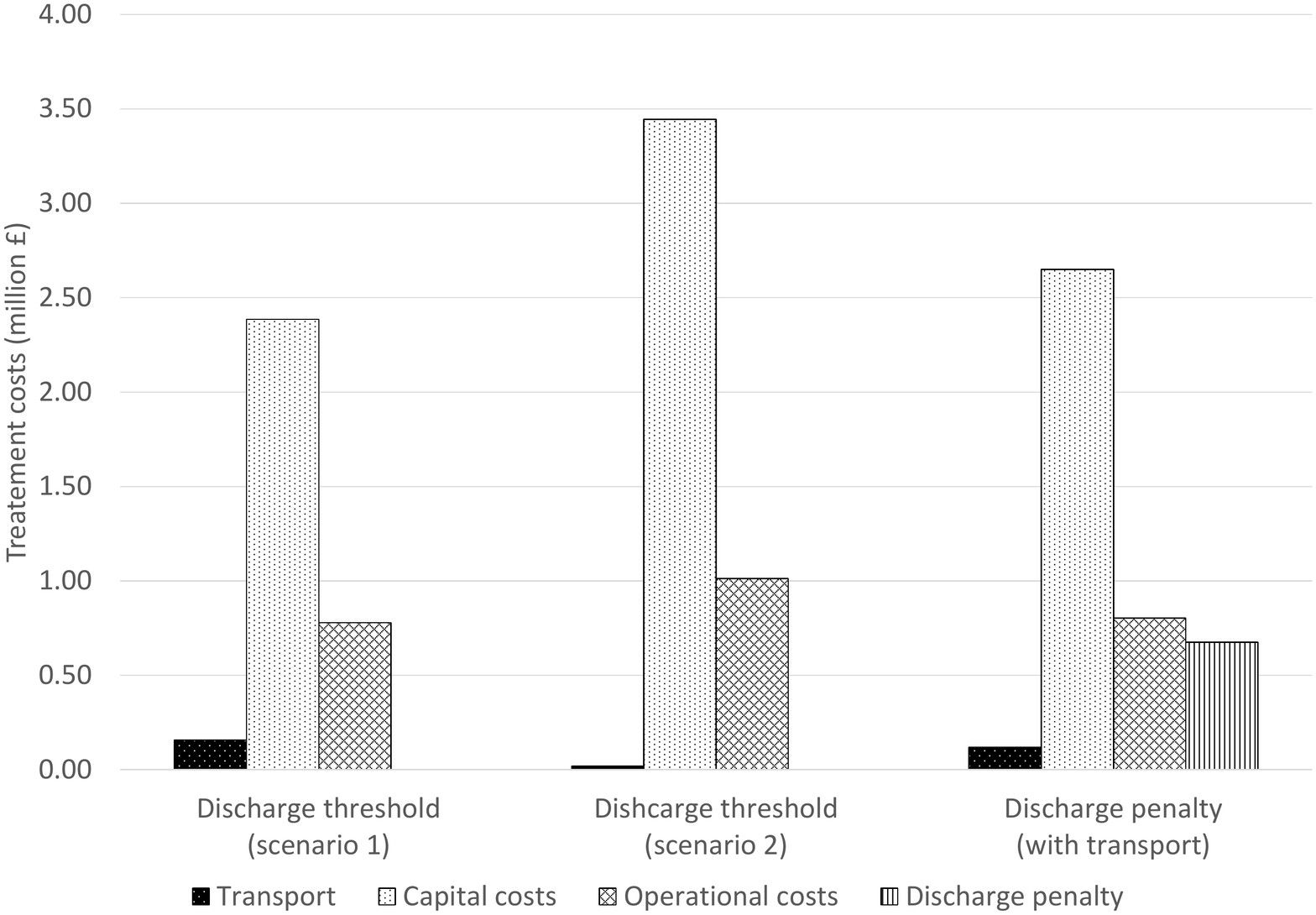}
    \caption{Decentralised treatment scheme layout}
    \label{fig:BarDischargeLimFin}
\end{figure}

The financial breakdown of these scenarios for construction and 10 years of operation is given in Fig.\ref{fig:BarDischargeLimFin}, including for comparison with the scenario from the previous section in which a discharge penalty is applied in place of hard constraints.

\begin{figure}[htbp]
    \centering
    \includegraphics[width=\mysize\textwidth]{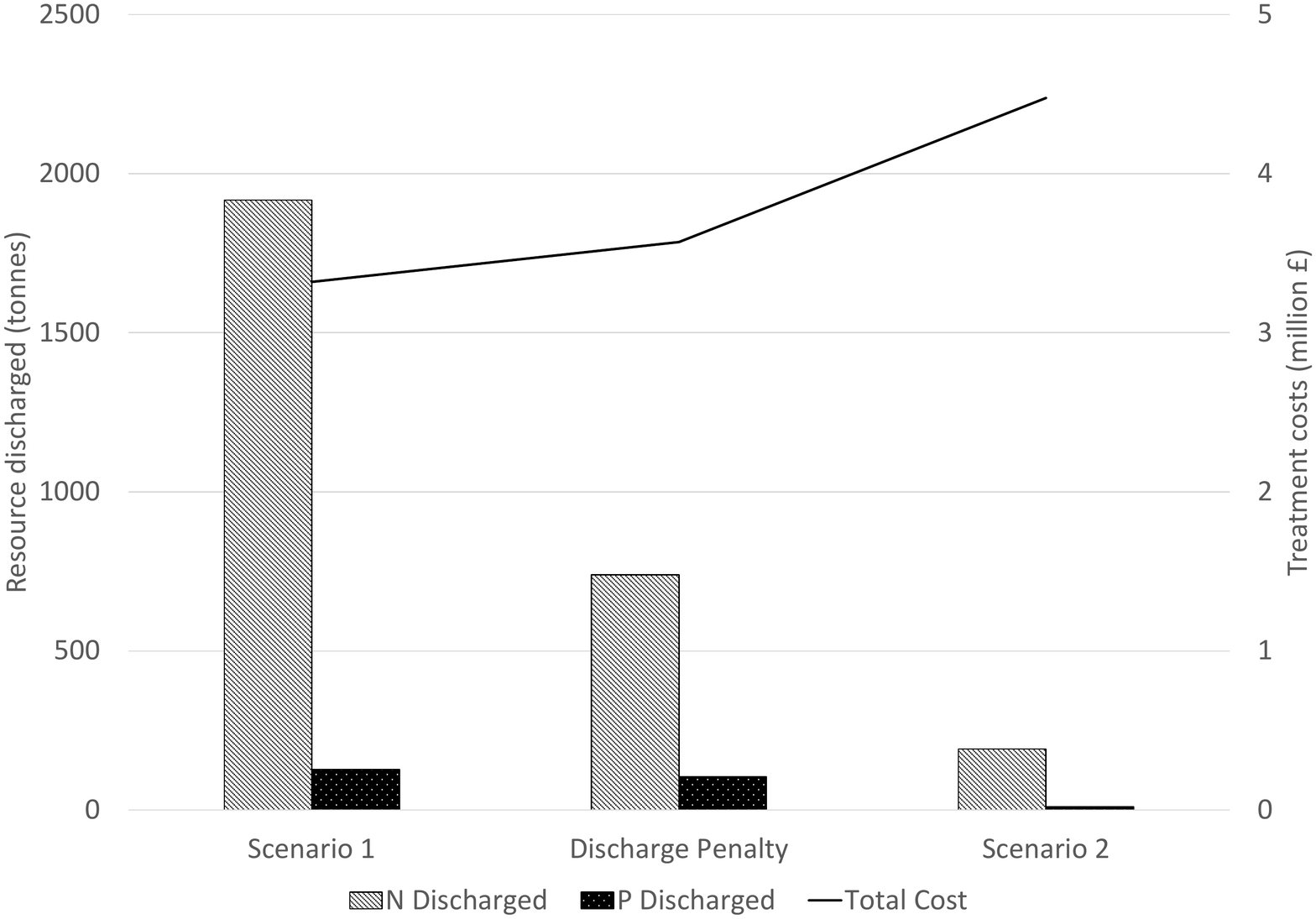}
    \caption{Decentralised treatment scheme layout}
    \label{fig:BarDischarge}
\end{figure}

As expected, to achieve the stricter limits on the environmental impact, an increased cost is incurred. This is illustrated in Fig.\ref{fig:BarDischarge}, in which the total quantity of P and N discharged is shown for the two scenarios as well scenario of Section \ref{Sec:Distributed} in which discharge penalties were applied. The total financial cost associated with each scenario is overlaid on the plot. While significant environmental improvements can be achieved by appropriate design, there is a financial barrier. The conflicting nature of these objectives underpines the need for a suitable design tool.

\subsection{Circular Economic Considerations}
With the benefits of resource-circular infrastructures well established \cite{EC_CircularEconomy}, eco-industrial park design should be carried out in a manner that best exploits the possibility of resource recovery and utilisation. Expanding on the previous examples, the production and sale of \(CH_{4}\) (as a biogas for the natural gas grid) and P and N (as a fertilizer) is now included in the formulation, in place of the discharge penalty used in previous scenarios. A selling price of \(\pounds0.16/m^3\) is assumed for methane generated as a gaseous biofuel \cite{ImperialCollegeLondon2018}. Similarly, selling prices of \(\pounds0.67/kg\) and \(\pounds0.27/kg\) are assumed for recovered N and P as nutrients respectively \cite{Alibaba2018}. The presence of a direct route to market for each of these is a large assumption for the design case, the financial implications of developing which lie outside the scope of this work. To represent a more conservative marketplace, a reduction factor of 0.5 is applied to the selling price.

\begin{figure}[htbp]
    \centering
    \includegraphics[width=\mysize\textwidth]{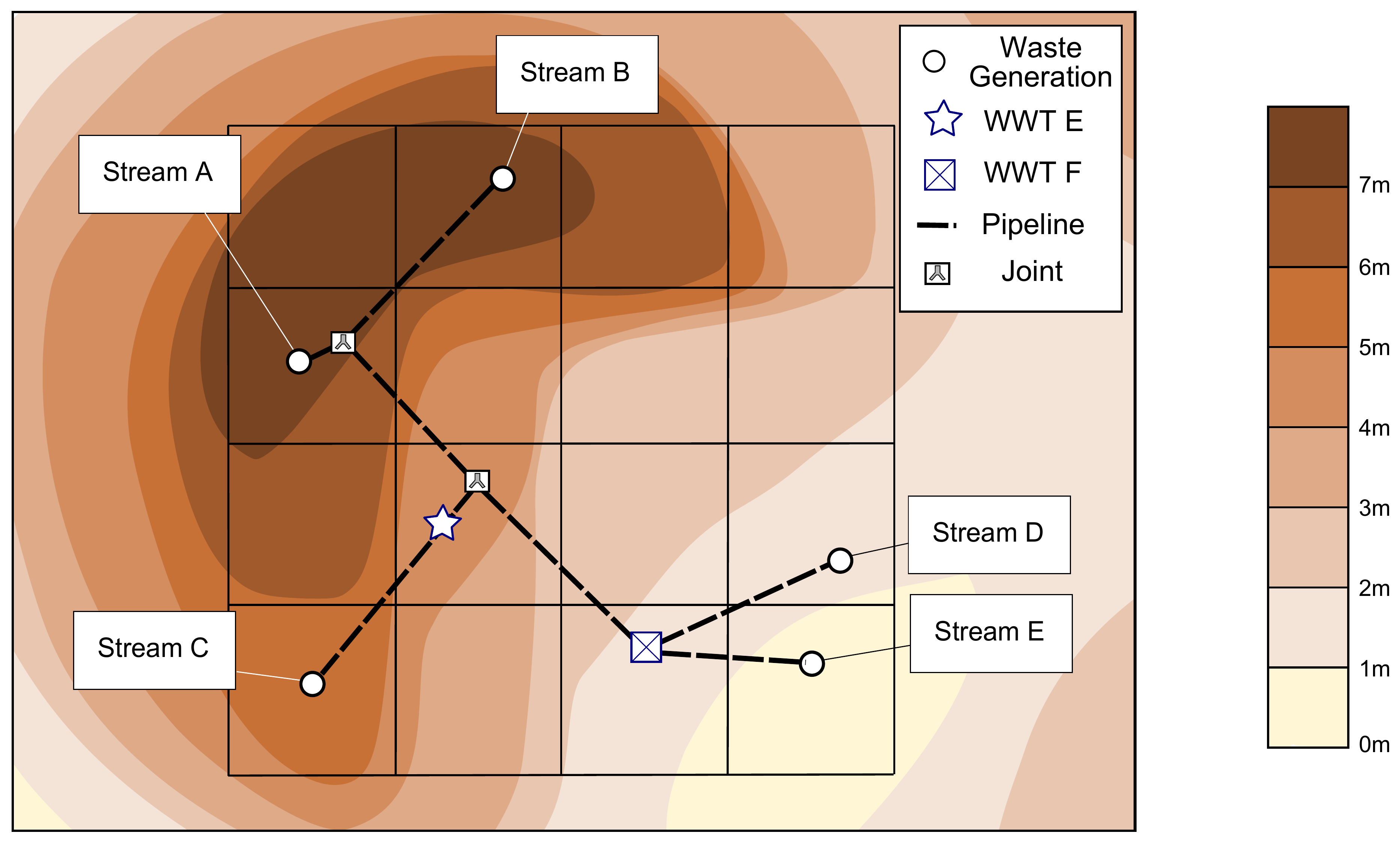}
    \caption{Design of treatment network with sale of recovered resources enabled}
    \label{fig:ContoursRecNoDisCost}
\end{figure}

The solution for this scenario is shown in Fig. \ref{fig:ContoursRecNoDisCost}. It can be seen that the P-recovery capability of the E-type plant now makes it the most financially attractive option, working in tandem here with an F-type plant. 

\begin{figure}[htbp]
    \centering
    \includegraphics[width=\mysize\textwidth]{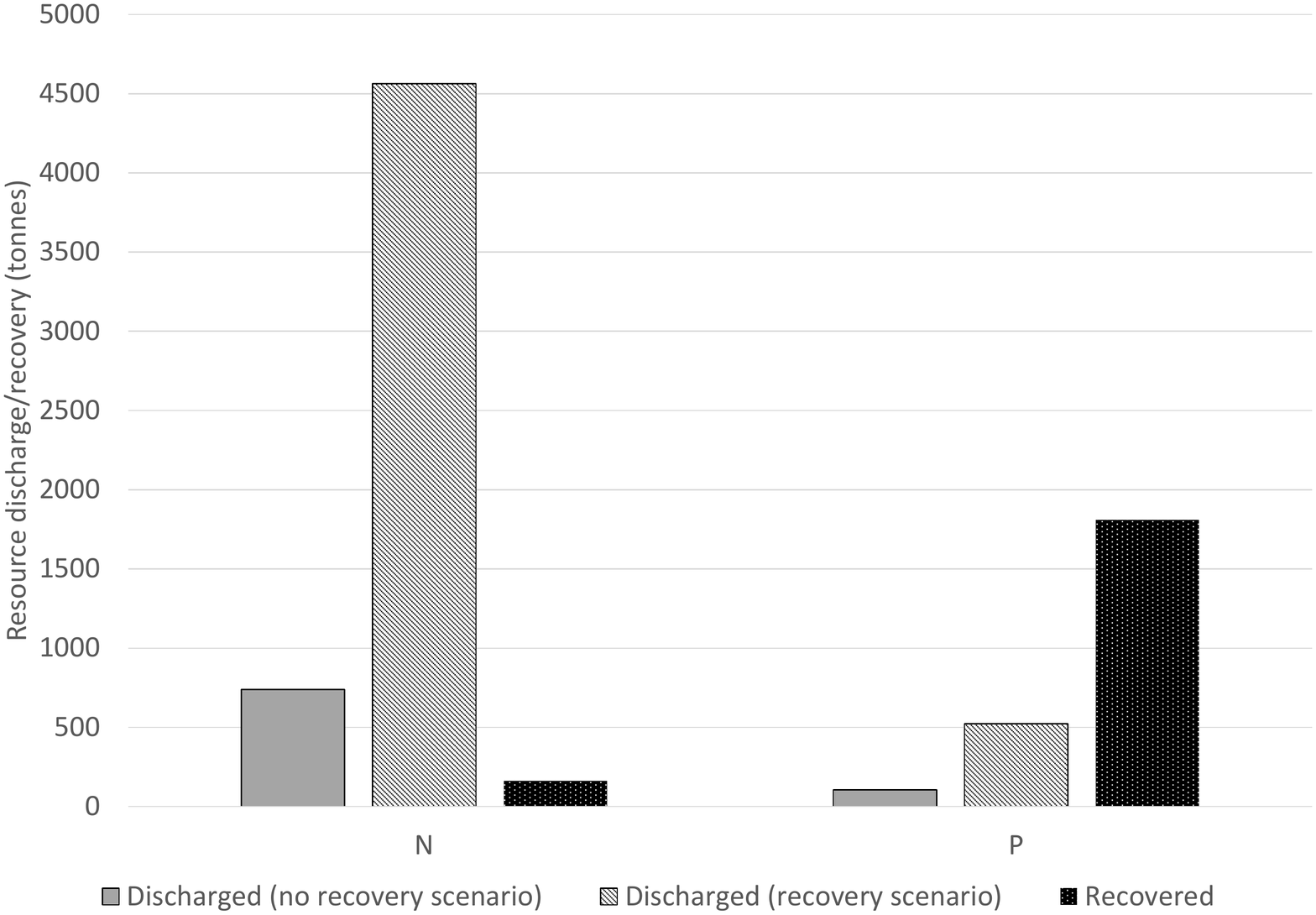}
    \caption{Discharged and recovered resources resulting from the design strategies of Fig. \ref{fig:ContoursDisCost} and Fig. \ref{fig:ContoursRecNoDisCost}}
    \label{fig:BarRecovery}
\end{figure}

Despite the introduction of resource recovery, a reduced environmental impact using this design cannot be assumed. The technology mix in this solution tends to result in a larger quantity of waste discharge to the environment, particularly due to the relatively low N-removal efficiency of the E-type plant. The total quantity of discharged and recovered N and P in 10 years of operation is shown in Fig. \ref{fig:BarRecovery}, with the discharge performance compared for reference to that of the scheme shown in Fig. \ref{fig:ContoursDisCost}. Although resources are now recovered, a clear increase in the quantity of discharged N and P can be seen.

\begin{figure}[htbp]
    \centering
    \includegraphics[width=\mysize\textwidth]{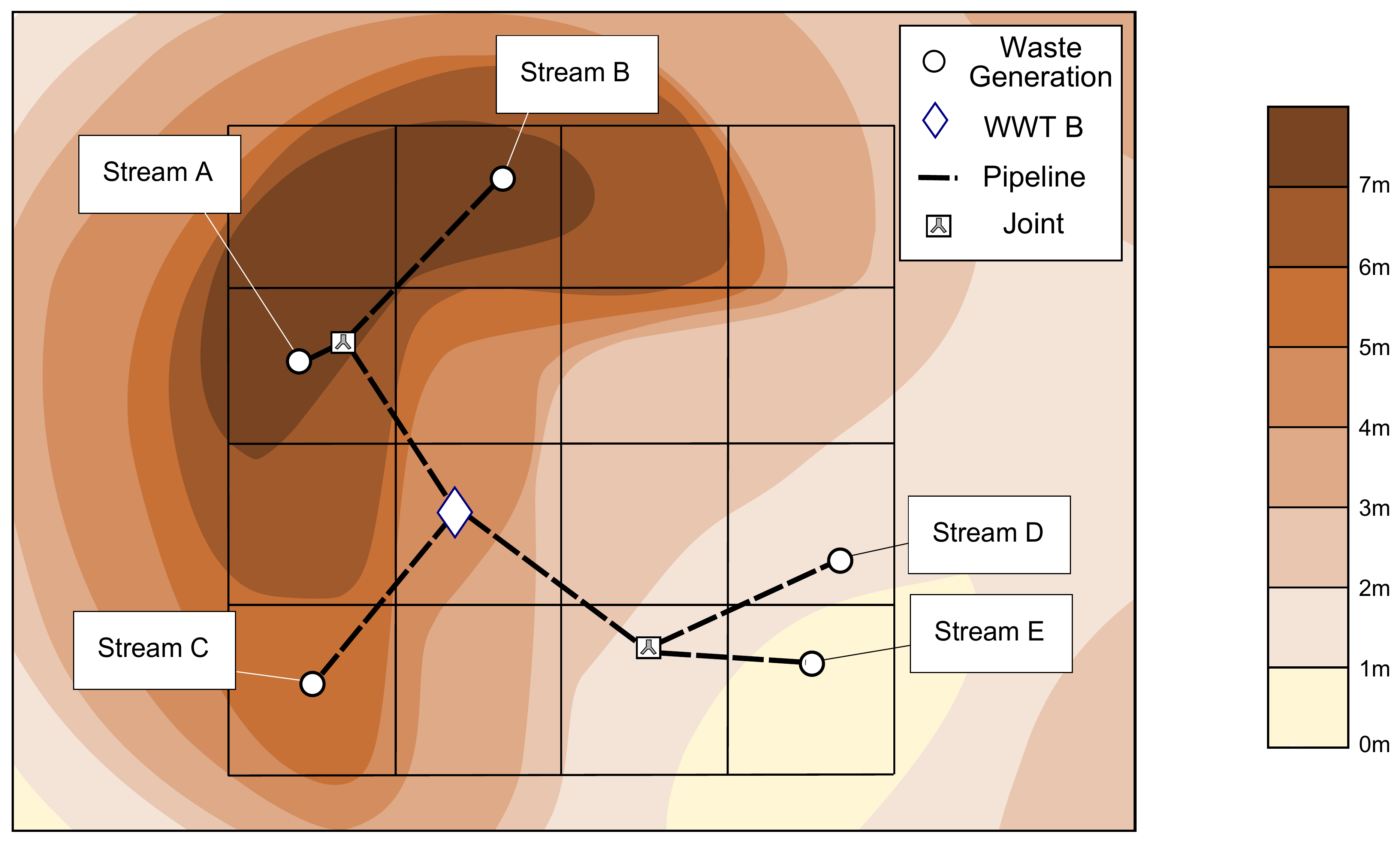}
    \caption{Design of treatment network with sale of recovered resources enabled and discharge penalties applied}
    \label{fig:ContoursRec}
\end{figure}

A new scenario is next presented in which the need to reduce the discharged waste while encouraging resource recovery is demonstrated by re-applying the discharge penalty while maintaining the resource recovery costs. The solution scheme is shown in Fig. \ref{fig:ContoursRec}.

Whereas in all other scenarios, small-scale treatment plants were favoured, in this case, it can be seen that all generated waste transported to a single centralised plant. The treatment B-type plant was selected which offers sufficient capacity to treat all streams, while deliver the highest \(CH_4\) recovery efficiency. The most financially appealing treatment solution in this scenario is to generate and sell biogas. This type of plant has one of the highest capital costs, however, by including resource recovery and appropriately penalising waste discharge, the environmental benefits of the technology are prioritised. The environmental performance of this scenario is summarised along with that of the scenario with a discharge penalty and no recovery and the scenario with only recovery and no discharge penalty in Table \ref{Tab:Case2}. These examples illustrate the multi-objective nature of the problem, and as such, the need for appropriate design tools to handle the competing requirements of the optimisation problem.

\begin{table}[htpb]
\centering
\begin{tabular}{l| c c c c c}
\hline
& \(N_{dis}\) & \(P_{dis}\) & \(N_{rec}\) & \(P_{rec}\) & \(CH_{4_{rec}}\) \\
 & \((tonnes)\) & \((tonnes)\) & \((tonnes)\) & \((tonnes)\) & \(\times10^{6}m^{3}\) \\
\hline\hline
Discharge penalty only & 740 & 105 & 0 & 0 & 0 \\
Recovery only & 4564 & 523 & 162 & 1810 & 101 \\
Recovery and discharge penalty & 3 & 0 & 0 & 0 & 116 \\
\hline
\end{tabular}
\caption{Cost of treatment facility over 10 operational years with recovered resources sold}
\label{Tab:Case2}
\end{table}

\section{Conclusion}\label{Conc}
Optimal design of wastewater treatment strategies in the context of environmentally sustainable eco-industrial parks is a complex challenge comprising of multiple design objectives across a spatially distributed mix of waste producers. The difficulty of the task is increased with the transition towards a more resource-circular treatment approach, in which resource recovery technologies are implemented to enable valuable by-products to be extracted. A generic optimisation formulation is proposed in this study to enable suitable design decisions to be made in this context. Key features of the formulation are the inclusion of transport pipeline network design in addition to technology selection. The formulation allows for the value of recovered resources to be considered, encouraging a shift towards more sustainable treatment design and operation. Furthermore, discharge penalties and strict limits to contaminant discharge can be applied at the park-level. A set of design case-studies was presented to illustrate the deployment of the optimisation formulation in different scenarios. The impact of different design assumptions on the preferred solution was emphasised, highlighting the importance of taking a whole-system perspective at the design stage. The case studies particularly implied that an optimal treatment strategy is highly sensitive to the value of the resources recovered and to the penalties applied to the discharge of different waste components. The optimisation formulation developed here allows for such factors to be considered, thereby promoting a transition towards more resource efficient wastewater treatment.

\section{Acknowledgements}\label{Ackn}
MG would like to acknowledge the UK Engineering and Physical Sciences Research Council (EPSRC) for providing financial support for funding our research via projects 'Resilient and Sustainable Biorenewable Systems Engineering Model' [EP/N034740/1].

\section*{References}\label{Ref}





\bibliographystyle{model1-num-names}
\bibliography{MainText_NoMarkup.bib}







\end{document}